\documentclass[11pt,a4paper]{article}
\usepackage[T1]{fontenc}
\usepackage[latin2]{inputenc}
\usepackage{ae,aecompl}
\usepackage{amsmath}
\usepackage{amsfonts}
\usepackage{amssymb}
\usepackage{amsthm}
\usepackage{graphicx}
\usepackage{color}
\usepackage[english,magyar]{babel}
\usepackage{listings}
\usepackage{enumerate}
\title{Finding weakly reversible realizations of chemical reaction networks using optimization}
\author{Gábor Szederkényi$^{a,1}$, Katalin M. Hangos$^{a,b}$, Zsolt Tuza$^{a,b}$}
\date{}
\newenvironment{gmatrix}[1]{\left[\begin{array}[c]{#1}}{\end{array}\right]}

\theoremstyle{definition}

\begin{document}
\selectlanguage{english} \maketitle
\begin{center}
$^{a}$
Computer and Automation Research Institute,\\
Hungarian Academy of Sciences\\
H-1518, P.O. Box 63, Budapest, Hungary\\
$^{b}$
Faculty of Information Technology, University of Pannonia\\
H-8200, Egyetem u. 10, Veszprém, Hungary
\end{center}

\medskip
\footnotetext[1]{The first author takes all responsibility for any errors or imperfections that might be found in this version of the document.}
\begin{abstract}
An algorithm is given in this paper for the computation of dynamically equivalent weakly reversible realizations with the maximal number of reactions, for chemical reaction networks (CRNs) with mass action kinetics. The original problem statement can be traced back at least 30 years ago \cite{Hars1981}. The algorithm uses standard linear and mixed integer linear programming, and it is based on elementary graph theory and important former results on the dense realizations of CRNs. The proposed method is also capable of determining if no dynamically equivalent weakly reversible structure exists for a given reaction network with a previously fixed complex set.
\end{abstract}

\newpage
\section{Introduction}
Chemical reaction networks (CRNs) are widely used for modeling chemical and biochemical processes in laboratory, industrial or natural environments \cite{Erdi1989}. Despite their simple mathematical structure, CRN models are capable of describing complex dynamical phenomena such as multiple steady states or oscillating behaviour that can be of key importance in understanding the operation of the modeled technological or living mechanisms \cite{Scott1994, Thomas2001, Son:2001}. From the point of view of mathematical sciences and systems theory, CRN models and the related analysis methods are also interesting and well-utilizable in the study of the characteristic features of nonlinear dynamical systems \cite{Angeli2009, Chellaboina2009}. A special class of CRNs is the set of networks whose dynamics is governed by the mass action law (MAL). In this paper, we will deal with such systems.

The foundations of chemical reaction network theory (CRNT) defining the basic notions and properties of reaction networks, were laid in the early papers and lecture notes \cite{Horn1972, Feinberg1977, Feinberg:79, Feinberg1987}. Since then, the development in the field has been continuous with significant contributions such as in \cite{Hars1981, Farkas1999, Siegel2000, c5, Craciun2005, Craciun2009, Craciun2010, Shinar2010}.

It has been known for long that CRNs with different structures and even with different set of chemical complexes may give rise to exactly the same set of differential equations describing the time-evolution of specie concentrations \cite{Horn1972}. This fact is called \textit{macro-equivalence} in \cite{Horn1972} and it is sometimes referred to as the "fundamental dogma of chemical kinetics". Different CRNs producing the same ODEs will be called \textit{dynamically equivalent} in this paper. Despite of the usefulness of dynamic equivalence in the analysis of CRNs, this property has only been studied to a rather limited extent in the literature. Conditions for the dynamical equivalence to a detailed balanced mechanism were first given in \cite{Averbukh1994}. In \cite{Craciun2008}, the authors gave necessary and sufficient conditions for the dynamic equivalence (also called \textit{confoundability}) of CRNs that was extended with a special case in \cite{Szederkenyi2009}. In \cite{Johnston2011}, the notion of dynamical equivalence was significantly extended by the introduction of \textit{conjugate networks} that mean CRNs with qualitatively the same dynamics where there is a well-defined mapping that takes trajectories of one system to into trajectories of the other. Also in \cite{Johnston2011}, conditions for conjugacy are given when the above mentioned mapping is a nontrivial linear transformation. 

In \cite{Szederkenyi2009b} the terms \textit{dense} and \textit{sparse realizations} were introduced for dynamically
equivalent reaction networks containing the maximal or minimal number of reactions with a fixed
set of complexes. Additionally, a computational method based on mixed integer linear
programming (MILP) was proposed to compute such realizations. In \cite{Szederkenyi2011}, important properties of dense realizations were described: it was shown that the structure (i.e. the unweighted directed
graph) of a dense realization of a given CRN is unique, and the unweighted directed reaction
graph of any other realization is the subgraph of the dense realization if the set of complexes is fixed. This is a key result that forms the basis of the algorithm presented in this paper. Furthermore, additional MILP based methods were also given in \cite{Szederkenyi2011} for computing reversible CRN realizations, and realizations with the minimal number of complexes. It was shown in \cite{Szederkenyi2010a} that the computation of dynamically equivalent detailed balanced and complex balanced CRNs can be traced back to simple linear programming (LP). 

Many strong results of Chemical Reaction Network Theory apply to weakly reversible networks. The most well-known examples are probably the \textit{Deficiency Zero} and \textit{Deficiency One Theorems} relating complex composition, network structure and the qualitative properties of CRN dynamics \cite{Feinberg1987}. 
Furthermore, it is known  that for any weakly reversible network structure, there always exists such a parametrization of the reaction rate coefficients that gives a complex balanced system \cite{Horn1972a}. Finally, the famous \textit{persistency conjecture}, says that the dynamics of weakly reversible networks is persistent in the sense that no trajectory that starts in the positive orthant has an $\omega$-limit point on the boundary of $\mathbb{R}^n_+$.

Similarly to numerous other key CRN properties such as deficiency, full reversibility, complex or detailed balance, weak reversibility is also not the property of the differential equations of the CRN models, but it is a realization property. Among the concluding open questions of \cite{Hars1981} (where the authors gave necessary and sufficient conditions for polynomial ODEs to be kinetic together with a constructive proof) we can read the following: "\textit{We may look for a mechanism in a class of mechanisms with a given - chemically relevant - property. Such a property may be conservativity, (weak) reversibility, zero deficiency or just structural stability as well.}" Additionally, in \cite{Johnston2011} the following is written: "\textit{The development of algorithms and computer software which can efficiently
check for viable weakly reversible target networks \dots is a primary interest}". Therefore it is of significant importance to work out a method that is able to construct a dynamically equivalent weakly reversible realization for a given CRN if it exists. 

Optimization tools based on e.g. linear/nonlinear programming or integer programming have been widely and successfully used in solving scientific and engineering problems \cite{Dantzig1997, Dantzig2003, Floudas1995}. Optimization methods can be very useful tools for deciding the feasibility of complex problems and give feasible solutions (if they exist), even if the original problem is  hard to treat analytically \cite{Boyd2004}. We will follow this approach in this paper.

The structure of the paper is the following. In section \ref{sec:CRNs} the applied models and relevant properties of CRNs will be summarized. Section \ref{sec:algorithm} contains the main contribution, i.e. the algorithm for determining weakly reversible realizations of CRNs. Illustrative examples are shown in section \ref{sec:examples}, while the conclusions can be found in section \ref{sec:conclusions}.

\section{Models and properties of chemical reaction networks}\label{sec:CRNs}
\subsection{Basic notions for the description of chemical reaction networks}
Following \cite{Feinberg:79} and several other works, we will characterize CRNs with the following three sets.
\begin{enumerate}
\item $\mathcal{S}=\{X_1, \dots, X_n \}$ is the set of \textit{species} or chemical substances. 
\item $\mathcal{C}=\{C_1,\dots,C_m \}$ is the set of \textit{complexes}. Formally, the complexes are represented as linear combinations of the species, i.e.
\begin{align}
C_i=\sum_{j=1}^n \alpha_{ij} X_j,~~~i=1,\dots,m
\end{align}
where $\alpha_{ij}$ are nonnegative integers and are called the \textit{stoichiometric coefficients}.
\item $\mathcal{R}=\{(C_i,C_j)~|~C_i,C_j\in\mathcal{C},~\text{and $C_i$ is transformed to $C_j$ in the CRN}\}$ is the set of \textit{reactions}. The relation $(C_i,C_j)\in\mathcal{R}$ will be denoted as $C_i\rightarrow C_j$. Moreover, a nonnegative weight, the \textit{reaction rate coefficient} denoted by $k_{ij}$ is assigned to each reaction $C_i\rightarrow C_j$.
\end{enumerate}  
The above characterization naturally gives rise to the following graph structure (often called 'Feinberg-Horn-Jackson graph' or simply reaction graph) of a reaction network \cite{Feinberg1987}. 
The weighted directed graph $G=(V,E)$ of a reaction network consists of a finite nonempty set $V$ of vertices and a finite set $E$ of ordered pairs of distinct vertices called directed edges. The vertices correspond to the complexes, i.e. $V=\{C_1,C_2,\dots C_m\}$, while the directed edges represent the reactions, i.e. $(C_i,C_j)\in E$ if complex $C_i$ is transformed to $C_j$ in the reaction network. The reaction rate coefficients $k_{ij}$  are assigned as positive weights to the corresponding directed edges $C_i \rightarrow C_j$ in the graph. 
A set of complexes $\{C_1,\dots,C_k \}$ is called a \textit{linkage class} of a CRN, if the complexes of the set are linked to each other in the reaction graph but not to any other complex. It is remarked that loops (i.e. directed edges that start and end at the same vertex) are not allowed in reaction graphs.
A reaction network is called \textit{reversible}, if whenever it contains the reaction $C_i\rightarrow C_j$, then the reverse reaction $C_j\leftarrow C_i$ is also present in the CRN. A reaction network is called \textit{weakly reversible}, if each complex in the reaction graph lies on at least one directed cycle (i.e. if complex $C_j$ is reachable from complex $C_i$ on a directed path in the reaction graph, then $C_i$ is reachable from $C_j$ on a directed path). 

A directed graph is called \textit{strongly connected} if there exists a path from each vertex of the graph to any other vertex. A \textit{strongly connected component} or simply \textit{strong component} of a directed graph is a set of vertices such that the directed edges between them provide a directed path from any vertex of the set to any other vertex, and to which no additional vertex can be added (i.e. a maximal strongly connected subgraph). Since a vertex is naturally reachable from itself through an empty path, strong components containing only one vertex will be called \textit{trivial strong components}. Clearly, the vertices of the individual linkage classes of a weakly reversible CRN form the strong components of the reaction graph.

Assuming mass-action kinetics, the following dynamical description will be used to describe the time-evolution of specie concentrations \cite{Feinberg:79, Feinberg1987}:
\begin{align}
\dot{x} & = Y\cdot A_k \cdot \psi(x)\label{eq:Feinberg_desc}
\end{align}
where  $x_i$ denotes the concentration of specie $X_i$. The $i$th column of $Y$ contains the composition of complex $C_i$, i.e. $[Y]_{ij}=\alpha_{ji}$. The structure and parameters of the reaction graph are stored in the column conservation matrix $A_k$ (also called the \textit{Kirchhoff matrix} of the CRN) as follows
\begin{equation}
[A_k]_{ij}=\left\{
\begin{array}{ccc}
-\sum_{l=1, l\ne i}^m k_{il}, & \text{if} & i=j \\
k_{ji}, & \text{if} & i\ne j
\end{array}
\right.
\end{equation}
Finally, $\psi:\mathbb{R}^n \mapsto\mathbb{R}^m$ is a monomial-type vector mapping defined by
\begin{equation} \label{Eq:monomials}
\psi_j(x)=\prod_{i=1}^n x_i^{[Y]_{ij}},~~~j=1,\dots,m
\end{equation}
\subsection{Dynamically equivalent realizations}
As it has been mentioned in the introduction, reaction networks with different structures and/or parametrizations can give rise to the same kinetic differential equations. Therefore, we will call two reaction networks given by the matrix pairs $(Y^{(1)}, A_k^{(1)})$ and $(Y^{(2)}, A_k^{(2)})$  \textit{dynamically equivalent}, if 
\begin{align}
Y^{(1)}A_k^{(1)}\psi^{(1)}(x)=Y^{(2)}A_k^{(2)}\psi^{(2)}(x)=f(x),~~\forall x\in\bar{\mathbb{R}}^n_+
\end{align}
where for $i=1,2$, $Y^{(i)}\in\mathbb{R}^{n\times m_i}$ have nonnegative integer entries, $A_k^{(i)}$ are valid Kirchhoff matrices, and
\begin{align}
\psi_j^{(i)}(x)=\prod_{k=1}^{m_i} x_k^{[Y^{(i)}]_{kj}},~~~i=1,2,~j=1,\dots,m.
\end{align}
In this case, $(Y^{(i)} A_k^{(i)})$ for $i=1,2$ are called \textit{realization}s of a kinetic vector field $f$ (see, e.g. \cite{Hars1981, Szederkenyi2010a} for more details). It is also appropriate to call $(Y^{(1)}, A_k^{(1)})$ a \textit{realization} of $(Y^{(2)}, A_k^{(2)})$ and vica versa.

\section{The algorithm for computing weakly reversible realizations}\label{sec:algorithm}
We recall that a \textit{dense realization} of a CRN contains the maximal number of reactions (i.e. maximal number of nonzero off-diagonal elements in $A_k$) with a given stoichiometric matrix $Y$ \cite{Szederkenyi2009b}. Based on the fact that with a given $Y$, all possible reactions are contained in the structurally unique dense realization \cite{Szederkenyi2011}, a straightforward idea is to try to find a dynamically equivalent weakly reversible 
mechanism starting from this superstructure. 

\subsection{Basic principle of the algorithm}\label{subsec:basicprinciple}
Very shortly, the underlying principle of the presented algorithm is that it only removes (if possible) from the dense realization
\begin{enumerate}[(i)]
\item edges that cannot be parts of any weakly reversible realization, 
\item edges the removal of which is necessarily implied by the deletion of edges belonging to set (i).
\end{enumerate}
The correct operation of the algorithm is based on the following two results.
\begin{itemize}
\item[\textbf{R1}] If each strongly connected component of a directed graph $G$ is contracted to a single vertex, the resulting directed graph is a directed acyclic graph \cite{Bang-Jensen2001}. (A directed graph is called \textit{acyclic} if it has no nontrivial strongly connected subgraphs.)
\item[\textbf{R2}] The structure of the dense realization of any CRN is unique, and the directed unweighted graph of any CRN realization is a subgraph of the directed unweighted graph of the dense realization, if the set of complexes is fixed \cite{Szederkenyi2011}.
\end{itemize}
From \textbf{R2} it follows that for obtaining a CRN superstructure including all possible structures, a dense realization must be computed. Directed edges between different strong components must be removed because they cannot lie on a directed cycle in any realization (\textbf{R1}). If this is not possible, then there is no weakly reversible realization of the CRN. However, if the deletion is possible, it may imply the removal of additional reactions because of the linear kinetic constraints (see eqs. \eqref{constr_1}-\eqref{constr_4}). In general, this may impair the weak reversibility of the obtained network.
In such a case, a new dense realization must be computed excluding the unnecessary edges identified in the previous step, and the procedure must be repeated until either a weakly reversible realization is found, or the deletion of undesired edges is no longer possible. In the latter case, no weakly reversible realization of the initial CRN exists with the given stoichiometric matrix $Y$.
\subsection{Definition of input data structure and the necessary additional procedures}
We will assume that an initial CRN realization is given with the matrices $(Y^{(0)}, A_k^{(0)})$, and $M=Y^{(0)}\cdot A_k^{(0)}$. The constraint set containing directed edges to be eliminated from the current realization is denoted as
\begin{align}
\mathcal{K}=\{(p_1,q_1),\dots,(p_s,q_s) \},~~~s<r.\label{eq:constraintK}
\end{align}
where $p_i$ and $q_i$ denote the indices of the initial and terminal vertices of the $i$th edge, respectively, and $r$ is the number of reactions in the CRN.

Now, the following simple procedure will be defined for later use.
\begin{align}
L=\text{\texttt{FindCrossComponentEdges}}(A_k^{in})
\end{align}
The input of the procedure is the Kirchhoff matrix $A_k^{in}$ of a CRN, and the output is a set $L$ containing the directed edges linking different strong components of the reaction graph. The strongly connected components of a directed graph can be determined in linear time using e.g. Kosaraju's, Tarjan's or Gabow's algorithm \cite{Bang-Jensen2001, Nuutila1994}. Moreover, the examined CRN is weakly reversible if and only if it contains at least 2 reactions, and the output $L$ of \texttt{FindCrossComponentEdges} is empty.

In the following part of this subsection, appropriate modifications of the algorithm described in \cite{Szederkenyi2009b} are given, adapted to the current problem.
\subsubsection{Constraints corresponding to mass action dynamics}
The main characteristics of mass action dynamics are taken into consideration in the form of the following constraints (see also \cite{Szederkenyi2009b} for more details) in both of the procedures presented in subsections \ref{subsub:remove} and \ref{subsubs:denseconstr}.
\begin{align}
& Y\cdot A_k = M \label{constr_1}\\
& \sum_{i=1}^m [A_k]_{ij}=0,~~~j=1,\dots,m \label{constr_2}\\
& [A_k]_{ij}\ge 0,~~i,j=1,\dots,m,~~i\ne j \\
& [A_k]_{ii}\le 0,~~i=1,\dots,m \label{constr_4}
\end{align}
where $Y$ and $M$ are known, and the decision variables are the off-diagonal elements of $A_k$. It's easy to see that constraints \eqref{constr_2}-\eqref{constr_4} represent the fact that we are searching for a valid Kirchhoff matrix.

\subsubsection{Checking whether a set of reactions is removable from a CRN realization}\label{subsub:remove}
We call a set of reactions $\mathcal{K}$ \textit{removable} from a CRN realization, 
 if there exists a dynamically equivalent CRN realization that does not contain the directed edges in $\mathcal{K}$.
To check this, it is worth separately defining the following LP-based and thus polynomial time procedure to avoid unnecessary MILP computations that are known to be NP-hard.

The constraints \eqref{constr_1}-\eqref{constr_4} are completed with the following ones
\begin{align}
[A_k]_{q_i,p_i}=0~~\text{for}~~i=1,\dots,s \label{constr:excl1}
\end{align}
where $(p_i,q_i)\in\mathcal{K}$ as it is written in eq. \eqref{eq:constraintK}.
The feasibility of the linear constraints \eqref{constr_1}-\eqref{constr_4} and \eqref{constr:excl1} can be checked by adding e.g. the following linear objective function to be minimized:
\begin{equation}
G_{LP}(A_k)=\sum_{\begin{array}{c} i,j=1 \\ i\ne j,~(j,i)\notin\mathcal{K} \end{array}}^m [A_k]_{ij}\label{eq:lin_obj}
\end{equation}
Clearly, eqs. \eqref{constr_1}-\eqref{eq:lin_obj} form a standard linear programming problem, the feasibility of which can be checked in polynomial time \cite{Dantzig1997}. 

Based on the above, we will define the procedure to check whether a set $\mathcal{K}$ of directed edges is removable from a CRN realization or not in the following way:
\begin{align}
F_{out} = \text{\texttt{IsRemovable}}(Y^{(i)}, A_k^{(i)}, \mathcal{K}), \label{eq:isremovable}
\end{align}
where the input data are as follows: $(Y^{(i)}, A_k^{(i)})$ is a CRN realization, and $\mathcal{K}$ is a constraint set of the form \eqref{eq:constraintK} containing the tail and head index pairs of the directed edges to be removed from $(Y^{(i)}, A_k^{(i)})$. The output $F_{out}$ is a boolean variable: its value is \texttt{true} if there exists a dynamically equivalent realization to $(Y^{(i)}, A_k^{(i)})$ not containing the edges listed in $K$, and it is \texttt{false} if there does not exist such realization.

\subsubsection{Computing the dense realization excluding given directed edges}\label{subsubs:denseconstr}
For the solution of this subtask, the well-known results on the combination of propositional logic with mixed integer programming are applied \cite{Raman1994, Bemporad1999}. According to these, a propositional logic problem, where a statement must be proved to be true if a set of compound statements are also given, can be solved through a linear integer problem.

If the procedure \texttt{IsRemovable} returns a \texttt{true} value, the dense realization of the CRN can be computed subject to the constraint that the edges listed in $\mathcal{K}$ are excluded from it. To define the corresponding MILP problem, first we add exactly the same linear constraints contained in eqs. \eqref{constr_1}-\eqref{constr:excl1} as in the previous case. To make the forthcoming problem computationally tractable, we also introduce the following bounds for the decision variables
\begin{align}
&[A_k]_{ij}\le u_{ij},~~u_{ij}>0,~~i,j=1,\dots,m, ~~i\ne j,~~(j,i)\notin \mathcal{K}\label{eq:bound1}\\
&[A_k]_{ii}\ge l_i,~~l_i<0,~~i=1,\dots,m\label{eq:bound2}
\end{align}
Here we are searching for such $A_k$ that contains the maximal number of nonzero off-diagonal elements. For this, logical variables denoted by $\delta$ are introduced and the following compound statements are constructed
\begin{align}
\delta_{ij}=1 \leftrightarrow [A_k]_{ij}>\epsilon,~~i,j=1,\dots,m,~~i\ne j,~~(j,i)\notin\mathcal{K} \label{compound01}
\end{align}
where the symbol "$\leftrightarrow$" denotes "if and only if", and $0<\epsilon\ll 1$ (i.e. elements of $A_k$ below $\epsilon$ are treated as zero). Considering also \eqref{eq:bound1}, statement \eqref{compound01} can be translated to the following linear inequalities 
\begin{align}
0\le [A_k]_{ij}-\epsilon \delta_{ij},~~i,j=1,\dots,m,~~i\ne j,~~(j,i)\notin\mathcal{K} \label{constr_ds_1}\\
0\le -[A_k]_{ij} + u_{ij}\delta_{ij},~~i,j=1,\dots,m,~~i\ne j~~(j,i)\notin\mathcal{K} \label{constr_ds_2}
\end{align}
Now it is possible to compute the realization containing the maximal number of reactions by maximizing the objective function
\begin{align}
G_{MILP}(\delta)=\sum_{\begin{small}\begin{array}{c}i,j=1\\i\ne j,~(j,i)\notin\mathcal{K}\end{array}\end{small}}^m \delta_{ij}.\label{ds_objfun}
\end{align}
Finally, the following procedure is defined based on the MILP problem given by eqs. \eqref{constr_1}-\eqref{constr:excl1}, \eqref{eq:bound1}-\eqref{eq:bound2}, and \eqref{constr_ds_1}-\eqref{ds_objfun}.
\begin{align}
A_k^d=\texttt{FindConstrDenseRealization}(Y^{(i)}, A_k^{(i)}, \mathcal{K}),
\end{align}
where the input data set is the same as in the case of the procedure \texttt{IsRemovable} (see eq. \eqref{eq:isremovable} and the corresponding description). The output $A_k^d$ is a dense realization that does not contain the directed edges listed in the set $\mathcal{K}$. If the set $\mathcal{K}$ is empty, then the procedure is the same that computes dense realizations and that was published in \cite{Szederkenyi2009b}.
\subsection{Formal description of the algorithm}\label{subsec:algdescr}
Now we are in the position to give the formal description of the procedure for determining weakly reversible CRN realizations. The input data of the procedure is an initial CRN realization $(Y^{(0)}$,$A_k^{(0)})$. The output is an $m\times m$ matrix that is the Kirchhoff matrix of the weakly reversible realization if there exists such, or a zero matrix if the procedure found no weakly reversible realizations. In the algorithm pseudocode, the auxiliary variable $ExitCondition$ is a boolean storing the exit condition from the main loop. The complete pseudocode of the procedure called \texttt{FindWeaklyReversibleRealization} with common notations and keywords can be found in Table \ref{tab:pseudocode}.

\begin{table}
\begin{center}
\framebox{
\parbox[c]{14cm}{
\begin{tabbing}
~~~~~ \= ~~~~\= ~~~~\= ~~~~\= ~~~~\= ~~~~ \= \kill
\underline{$A_k^{out}$=\texttt{FindWeaklyReversibleRealization}($Y^{(0)}$,$A_k^{(0)}$)}\\
1 \> $A_k^{out}$:=$0\in\mathbb{R}^{m\times m}$; $ExitCondition$:=\texttt{false};\\
2 \> $Y:=Y^{(0)}$; $A_k:=A_k^{(0)}$; $F_{out}$:=\texttt{true}; $\mathcal{K}:=\{\}$; $L:=\{\}$; \\
3 \> \texttt{while} ($ExitCondition$=\texttt{false}) \texttt{do}\\
4 \> \texttt{begin}\\
5 \> \> \texttt{if} ($\mathcal{K}\ne \{\}$) \texttt{then} $F_{out}$:=\texttt{IsRemovable}($Y$,$A_k$,$\mathcal{K}$);\\
6 \> \> \texttt{if} ($F_{out}=$\texttt{true}) \texttt{then}\\
7 \> \> \texttt{begin}\\
8 \> \> \> $A_k$:=\texttt{FindConstrDenseRealization}($Y$,$A_k$,$\mathcal{K}$);\\
9 \> \> \> $L$:=\texttt{FindCrossComponentEdges}($A_k$);\\
10 \> \> \> \texttt{if} ($L=\{\}$) \texttt{then} $ExitCondition$:=\texttt{true}; $A_k^{out}$:=$A_k$;\\
11 \> \> \> \texttt{else} $\mathcal{K}:=\mathcal{K}\cup L$;\\
12 \> \> \texttt{end} \\
13 \> \> \texttt{else} $ExitCondition$:=\texttt{true};\\
14 \> \texttt{end}\\
15 \> \texttt{return} $A_k^{out}$;
\end{tabbing}
}}
\end{center}
\caption{Pseudocode of the algorithm for finding weakly reversible realizations}\label{tab:pseudocode}
\end{table}

\subsection{Main properties of the algorithm}
The above described algorithm always finds a dynamically equivalent weakly reversible realization, if it exists. This clearly follows from the results \textbf{R1}, \textbf{R2} and the basic principles of the algorithm described in subsection \ref{subsec:basicprinciple}. From these facts it also follows that the algorithm finds the 'densest' weakly reversible realization that structurally contains any other weakly reversible realizations (for illustration, see the results in subsection \ref{subsec:example1}). 
An apparent drawback of the algorithm that it contains a step where a MILP problem has to be solved that is known to be NP-hard. However, the parallel (columnwise) implementation of the procedure \texttt{FindConstrDenseRealization} is possible as it is explained in \cite{Szederkenyi2009b}, that significantly extends the number of complexes that can be handled in practice. We remark that the immediate deletion of the columns/rows corresponding to the isolated complexes from matrices $Y$ and $A_k$ after calling the procedure \texttt{FindConstrDenseRealization} is also possible, but this requires the renumbering of complexes. This is a technical detail of implementation and does not affect the principle or final output of the algorithm. (However, decreasing the number of optimization variables and constraints in each step in such a way might be required especially in the case of larger networks.)
\section{Examples}\label{sec:examples}
The following examples were implemented in the MATLAB R13 computation environment using the YALMIP and the Multi-Parametric Toolboxes \cite{Loefberg2004, Kvasnica2004}. The examples were run on a desktop PC with dual Intel Xeon 1.8GHz CPU and 2 Gigabytes of RAM. The implementation of dense realization computation was not parallel, therefore all the decision variables and constraints were put into a single optimization problem in procedure \texttt{FindConstrDenseRealization}. The strong components of the reaction graphs were identified using Kosaraju's algorithm \cite{Sharir1981}.
\subsection{Weakly reversible realizations of a simple irreversible network}\label{subsec:example1}
\begin{figure}[!htbp]%
\centering
\includegraphics[width=12cm]{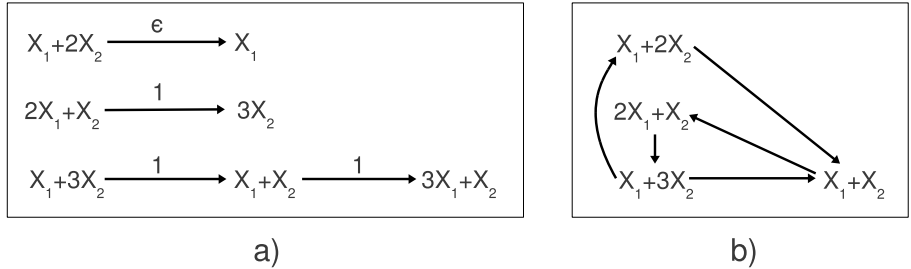}
\caption{a) Simple irreversible network from \cite{Johnston2011}, b) Structure of one of its possible weakly reversible realizations determined in \cite{Johnston2011}}
\label{JohnstonSiegel_00}
\end{figure}
\begin{figure}[!htbp]%
\centering
\framebox{\includegraphics[width=9cm]{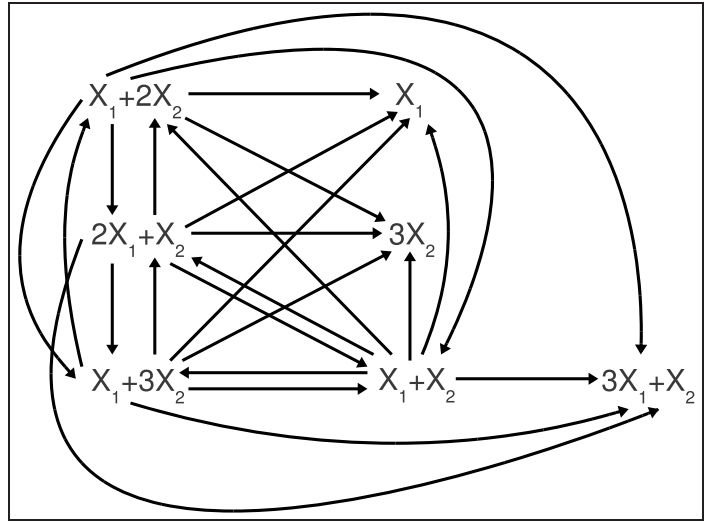}}
\caption{Structure of the dense realization of the reaction network shown in Fig. \ref{JohnstonSiegel_00} a)}
\label{JohnstonSiegel_01}
\end{figure}
\begin{figure}[!htbp]%
\centering
\framebox{\includegraphics[width=9cm]{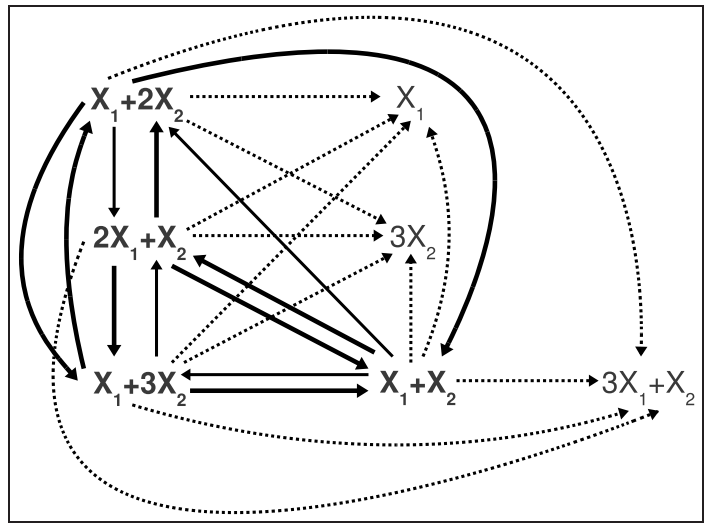}}
\caption{Illustration of the operation of the algorithm on the example given in section \ref{subsec:example1}}
\label{JohnstonSiegel_02}
\end{figure}
\begin{figure}[!htbp]%
\centering
\framebox{\includegraphics[width=7cm]{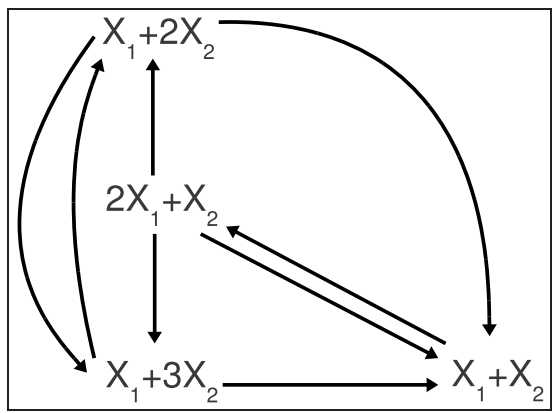}}
\caption{The structure of the obtained dynamically equivalent weakly reversible CRN}
\label{JohnstonSiegel_03}
\end{figure}

The simple network that can be seen in Fig. \ref{JohnstonSiegel_00} a) was taken from \cite{Johnston2011} (Example 3). In \cite{Johnston2011}, it is shown that for any positive $\epsilon$, the network has a possible weakly reversible realization with the structure shown in Fig. \ref{JohnstonSiegel_00}b). (The computation and analysis of the parameters of the CRN in Fig. \ref{JohnstonSiegel_00}b can be found in \cite{Johnston2011}.) The stoichiometric matrix of the network is
\begin{equation}
Y=\begin{gmatrix}{ccccccc}
1 & 1 & 2 & 0 & 1 & 1 & 3 \\
2 & 0 & 1 & 3 & 3 & 1 & 1  
\end{gmatrix}
\end{equation}
The nonzero elements of the Kirchhoff matrix $A_k\in\mathbb{R}^{7\times 7}$ with $\epsilon=1.5$ are
\begin{align}
[A_k]_{2,1}=1.5,~[A_k]_{4,3}=1,~[A_k]_{6,5}=1,~[A_k]_{7,6}=1~.
\end{align}
For this network, the algorithm described in section \ref{subsec:algdescr} and Table \ref{tab:pseudocode} works as follows. After the initialization steps, the dense realization containing all possible reactions with an empty constraint set $\mathcal{K}$ is computed (line 8 of the pseudocode). The Kirchhoff matrix of the dense realization is given by
\begin{align}
A_k^{(1)} =
\begin{gmatrix}{ccccccc}
-1.25   &      0  &  0.1   &      0  &  0.1  &  0.1 & 0 \\
    0.55   &      0  &  0.1  &       0  &  0.4333  &  0.5 & 0 \\
    0.1   &      0  & -1.4   &      0  &  0.1  &  0.1 & 0 \\
    0.3    &     0  &  0.8    &     0  &  0.3  &  0.1 & 0 \\
    0.1    &     0  &  0.2     &    0  & -1.1333  &  0.1 & 0 \\
    0.1    &     0  &  0.1     &    0  &  0.1  & -1.9 & 0 \\
    0.1    &     0  &  0.1     &    0  &  0.1  &  1 & 0 \\
\end{gmatrix}
\end{align}
The structure of the dense realization is shown in Fig. \ref{JohnstonSiegel_01}.
The following steps can be followed using Fig. \ref{JohnstonSiegel_02}. The dense realization is not weakly reversible because there are edges between different strong components (line 9). The complexes of the single nontrivial strong component are indicated by boldface labels in Fig. \ref{JohnstonSiegel_02}. It is clear from the figure that the edges adjacent to the complexes $X_1$, $3X_2$, $3X_1 + X_2$ are to be removed. These edges are drawn 
with dotted arrows in the figure. Thus, the constraint list is 
\[
\mathcal{K}=\{(1,2),(1,4),(1,7),(3,2),(3,4),(3,7),(5,2),(5,4),(5,7),(6,2),(6,4),(6,7) \}.
\] 
The next iteration of the algorithm (line 5) gives that these edges can be removed from the realization. Now a new dense realization is computed excluding edges in $\mathcal{K}$ (line 8). The reactions of this dense realization are indicated by thick arrows in Fig. \ref{JohnstonSiegel_02}. Note that the constraint of excluding the reactions in $\mathcal{K}$ from the network resulted in the removal of directed edges $X_1+X_2\rightarrow X_1+2X_2$, $X_1+2X_2\rightarrow 2X_1+X_2$, $X_1+3X_2\rightarrow 2X_1+X_2$,  and $X_1+X_2\rightarrow X_1+3X_2$, too, that were within the nontrivial strong component of the previous step. In this case, the resulting network remained weakly reversible (lines 9-10), so the algorithm can stop with success and return the determined dynamically equivalent weakly reversible CRN. The structure of the resulting  network is shown in Fig. 
\ref{JohnstonSiegel_03}. The Kirchhoff matrix of the obtained realization is the following
\begin{align}
A_k^{(2)} =
\begin{gmatrix}{ccccccc}
 -3.2   &      0  &   1.8    &     0  &  0.1    &     0 & 0\\
         0    &     0    &     0    &     0    &     0    &     0 & 0\\
         0   &      0  & -2  &       0    &     0  &  2 & 0 \\
         0   &      0   &      0    &     0    &     0    &     0 & 0 \\
    0.1    &     0  &  0.1   &      0  & -1.05   &      0 & 0 \\
    3.1   &      0  &  0.1   &      0  &  0.95 &  -2 & 0 \\
         0   &      0   &      0    &     0   &      0   &      0 & 0
\end{gmatrix}
\end{align}
The running time of the algorithm was 11s using the hardware/software environment described in the beginning of section \ref{sec:examples}.

\begin{figure}[!htbp]%
\centering
\framebox{\includegraphics[width=6cm]{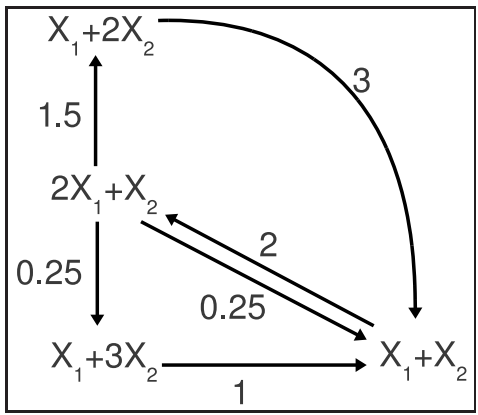}}
\caption{Complex balanced realization of the CRN described in section \ref{subsec:example1}}
\label{JohnstonSiegel_04}
\end{figure}

It is interesting to note that the obtained weakly reversible realization is not complex balanced. (For convenience, we briefly define the notion of complex balance:
A CRN realization $(Y,A_k)$ is called \textit{complex balanced} if there exists a positive vector  $x^*\in\mathbb{R}^n_+$ for which $A_k\psi(x^*)=0$ \cite{Horn1972, Siegel2000, Craciun2009}.)  
However, using the polynomial-time algorithm described in \cite{Szederkenyi2010a}, we can compute a complex balanced realization with 6 reactions in 0.1s, that is shown in Fig. \ref{JohnstonSiegel_04}. It is easy to see that the unweighted directed graphs of Figs. \ref{JohnstonSiegel_00} b) and Fig. \ref{JohnstonSiegel_04} are the proper subgraphs of the structure visible in Fig. \ref{JohnstonSiegel_03}.
\subsection{Weakly reversible realization of a kinetic polynomial system}
\begin{figure}[!htbp]%
\centering
\framebox{\includegraphics[width=6cm]{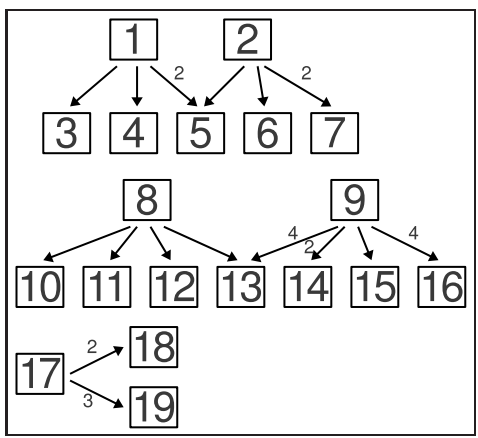}}
\caption{CRN realizing the dynamics of eq. \eqref{eq:example2}. Only reaction rates different from 1 are indicated.}
\label{react_02_original}
\end{figure}

The equations and initial CRN in this example were also used in \cite{Szederkenyi2010a} and \cite{Szederkenyi2010b} for illustrating other computation methods. We start from the polynomial ordinary differential equations given by
\begin{align}
\dot{x}_1 & = x_3^2 - x_1 x_2 + x_3 x_4 - 2 x_1 x_2^2 x_3 \nonumber \\
\dot{x}_2 & = x_3^2 - x_1 x_2 + 2 x_3 x_4 - 4 x_1 x_2^2 x_3 \nonumber \\
\dot{x}_3 & = -2 x_3^2 + x_1 x_2 - x_1 x_2^2 x_3 + 2 x_4^3 \label{eq:example2}\\
\dot{x}_4 & = x_1 x_2 - x_3 x_4 + 4 x_1 x_2^2 x_3 - 3 x_4^3 \nonumber
\end{align}
Using the inverse kinetic algorithm first published in \cite{Hars1981} (see also \cite{Szederkenyi2010a} for a summary about the kinetic realizability of polynomial vector fields), a CRN shown in Fig. \ref{react_02_original} can be constructed algorithmically that gives the dynamics \eqref{eq:example2}. 
The numbering of complexes in the figure is the following:
\begin{align}
& 1: 2  X_3,~
2:  X_3+  X_4,~
3:  X_1+2  X_3,~
4:  X_2 + 2  X_3,\nonumber\\
& 5:  X_3,~
6:  X_1 +   X_3 +   X_4,~
7:  X_2 +   X_3 +   X_4,~\nonumber\\
& 8:  X_1+  X_2,~
9:  X_1 + 2  X_2 +   X_3,~
10:  X_1,~
11:  X_2,~\nonumber\\
& 12:  X_1 +   X_2 +   X_4,~
13:  X_1 +   X_2 +   X_3,~
14:2  X_2+  X_3,~\nonumber\\
& 15:  X_1+2  X_2,~
16:   X_1 + 2  X_2 +   X_3 +   X_4,~\nonumber\\
& 17: 3  X_4,~18:   X_3+3   X_4,~ 19: 2   X_4 
\end{align}

The dense realization of the network contains 80 reactions, therefore it is not shown in a figure. For the sake of completeness, the list of the reactions (i.e. weighted directed edges) of the dense realization in the form (\textit{source vertex number, destination vertex number, rate coefficient}) is given below:\\
(5, 1, 0.5), (7, 1, 0.1), (11, 1, 0.1), (15, 1, 0.8), (1, 2, 0.1), (3, 2, 0.1), (5, 2, 0.1), (7, 2, 0.1), (12, 2, 0.3), (1, 3, 0.1), (5, 3, 0.1), (7, 3, 0.1), (12, 3, 0.1), (1, 4, 0.1), (3, 4, 0.4), (5, 4, 0.1), (7, 4, 0.1), (11, 4, 0.1), (3, 5, 0.1), (7, 5, 0.1), (11, 5, 0.1), (15, 5, 0.1), (3, 6, 0.1), (5, 6, 0.1), (7, 6, 0.1), (1, 7, 0.1), (3, 7, 0.1), (5, 7, 0.2), (12, 7, 0.3), (1, 8, 0.1), (3, 8, 0.1), (5, 8, 0.3), (7, 8, 0.3), (11, 8, 1.2), (1, 9, 0.1), (5, 9, 0.1), (7, 9, 0.1), (11, 9, 0.2), (1, 10, 0.5), (3, 10, 0.8), (5, 10, 0.1), (7, 10, 0.1), (12, 10, 0.1), (3, 11, 0.1), (5, 11, 0.1), (7, 11, 0.1), (1, 12, 0.1), (3, 12, 0.1), (5, 12, 0.1), (7, 12, 0.1), (1, 13, 0.1), (3, 13, 0.1), (5, 13, 0.1), (7, 13, 0.1), (11, 13, 0.1), (15, 13, 0.1), (1, 14, 0.1), (3, 14, 0.1), (5, 14, 0.1), (7, 14, 0.1), (12, 14, 0.1), (3, 15, 0.1), (5, 15, 0.1), (7, 15, 0.7), (11, 15, 0.1), (5, 16, 0.25), (7, 16, 0.3), (11, 16, 0.7), (15, 16, 0.2), (3, 17, 0.1), (5, 17, 0.1), (7, 17, 0.1), (3, 18, 0.1), (5, 18, 0.1), (7, 18, 0.4), (3, 19, 0.1), (5, 19, 0.1), (7, 19, 0.1), (11, 19, 0.1), (15, 19, 0.1).

The dense realization has 13 strong components. Out of these, there is only one nontrivial strongly connected component containing the vertices 1, 3, 5, 7, 11, 12, and 15. After identifying all directed edges linking different strong components, we obtain that the following 53 edges given in the form (\textit{source vertex number, destination vertex number}) should be deleted from the dense realization in the next step of the algorithm:\\
(1,2), (3,2), (5,2), (7,2), (12,2), (1,4), (3,4), (5,4), (7,4), (11,4), (3,6), (5,6), (7,6), (1,8), (3,8), (5,8), (7,8), (11,8), (1,9), (5,9), (7,9), (11,9), (1,10), (3,10), (5,10), (7,10), (12,10), (1,13), (3,13), (5,13), (7,13), (11,13), (15,13), (1,14), (3,14), (5,14), (7,14), (12,14), (5,16), (7,16), (11,16), (15,16), (3,17), (5,17), (7,17), (3,18), (5,18), (7,18), (3,19), (5,19), (7,19), (11,19), (15,19).

\begin{figure}[!htbp]%
\centering
\framebox{\includegraphics[width=6cm]{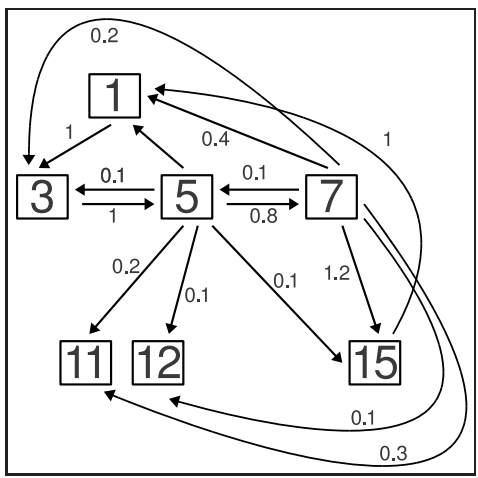}}
\caption{CRN realization of eq. \eqref{eq:example2} containing the vertices of the only nontrivial strong component of the dense realization}
\label{react_02_step1}
\end{figure}

\begin{figure}[!htbp]%
\centering
\framebox{\includegraphics[width=6cm]{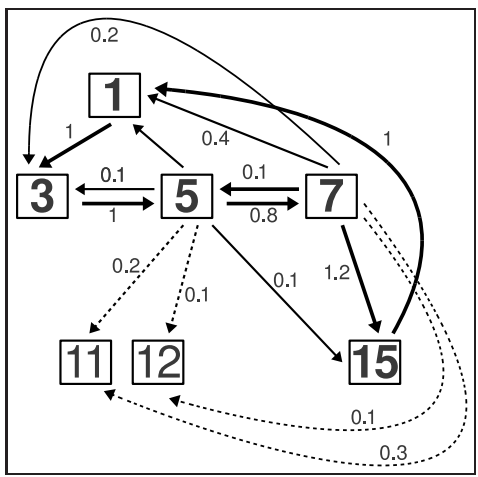}}
\caption{Illustration of the final step of the algorithm on the CRN corresponding to eq. \eqref{eq:example2}}
\label{react_02_step1_1}
\end{figure}

\begin{figure}[!htbp]%
\centering
\framebox{\includegraphics[width=6cm]{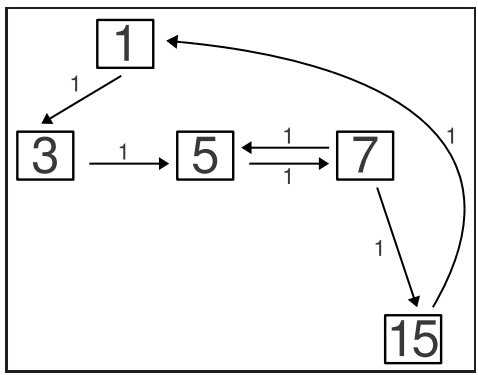}}
\caption{The obtained weakly reversible realization of eq. \eqref{eq:example2}}
\label{react_02_step2}
\end{figure}

All the above listed edges were possible to remove from the dense realization (their removal implied the deletion of 12 additonal reactions), and the resulting constrained dense realization is shown in Fig.\ref{react_02_step1}. It is clear from the figure that edges adjacent to vertices no. 11 and 12 have to be removed in the following step. This final step is illustrated in Fig. \ref{react_02_step1_1}, where the meaning of the line types is the same as in the case of Fig. \ref{JohnstonSiegel_02}. With the removal of edges (5,11), (5,12), (7,11), (7,12), the directed edges (5,1), (7,1), (7,3), (5,3), (7,5) and (5,15) were also deleted, and the resulting weakly reversible realization (of deficiency 0) is shown in Fig. \ref{react_02_step2}. The total running time of the algorithm was 80.5s.
\section{Conclusions}\label{sec:conclusions}
A numerical algorithm for finding dynamically equivalent we\-akly reversible realizations of chemical reaction networks was proposed in this paper. The algorithm computes a maximal superstructure that contains all other possible weakly reversible structures as proper subgraphs, and it is able to determine if no weakly reversible realization exists. The method uses linear and mixed integer linear programming steps and it is based primarily on the computation and properties of dense realizations published in \cite{Szederkenyi2009b, Szederkenyi2011}. The computationally critical MILP step of the algorithm can be implemented parallelly. The operation of the algorithm has been illustrated on numerical examples.
\section*{Acknowledgements}
This research work has been partially supported by the Hungarian Scientific Research Fund
 through grant no. K67625 and by the Control Engineering Research Group of the Budapest
 University of Technology and Economics. 

\bibliography{react_weakly}
\end{document}